\documentclass[12pt,reqno]{amsart}
\usepackage{amssymb}
\title[Numerical bounds for semi-stable families]
{Numerical bounds for semi-stable families of curves or of certain higher dimensional manifolds}
\author[Eckart Viehweg]{Eckart Viehweg}
\address{Universit\"at Duisburg-Essen, Mathematik, 45117 Essen, Germany}
\email{viehweg@uni-essen.de}
\thanks{This work has been supported by the ``DFG-Schwerpunktprogramm
Globale Methoden in der Komplexen Geometrie'', and by the DFG-Leibniz program.}
\author[Kang Zuo]{Kang Zuo}
\address{Universit\"at Mainz,
Fachbereich 17, Mathematik,
55099 Mainz, Germany}
\email{kzuo@mathematik.uni-mainz.de}
\begin{document}
\theoremstyle{plain}
\newtheorem{thm}{Theorem}[section]
\newtheorem{theorem}[thm]{Theorem}
\newtheorem{lemma}[thm]{Lemma}
\newtheorem{corollary}[thm]{Corollary}
\newtheorem{proposition}[thm]{Proposition}
\newtheorem{construction}[thm]{Construction}
\newtheorem{addendum}[thm]{Addendum}
\theoremstyle{definition}
\newtheorem{notations}[thm]{Notations}
\newtheorem{notation}[thm]{Notation}
\newtheorem{problem}[thm]{Problem}
\newtheorem{remark}[thm]{Remark}
\newtheorem{remarks}[thm]{Remarks}
\newtheorem{definition}[thm]{Definition}
\newtheorem{claim}[thm]{Claim}
\newtheorem{assumption}[thm]{Assumption}
\newtheorem{setup}[thm]{Set-up}
\newtheorem{assumptions}[thm]{Assumptions}
\newtheorem{properties}[thm]{Properties}
\newtheorem{example}[thm]{Example}
\newtheorem{examples}[thm]{Examples}
\newtheorem{conjecture}[thm]{Conjecture}
\numberwithin{equation}{section}
\catcode`\@=11
\def\opn#1#2{\def#1{\mathop{\kern0pt\fam0#2}\nolimits}}
\def\bold#1{{\bf #1}}%
\def\underrightarrow{\mathpalette\underrightarrow@}
\def\underrightarrow@#1#2{\vtop{\ialign{$##$\cr
 \hfil#1#2\hfil\cr\noalign{\nointerlineskip}%
 #1{-}\mkern-6mu\cleaders\hbox{$#1\mkern-2mu{-}\mkern-2mu$}\hfill
 \mkern-6mu{\to}\cr}}}
\let\underarrow\underrightarrow
\def\underleftarrow{\mathpalette\underleftarrow@}
\def\underleftarrow@#1#2{\vtop{\ialign{$##$\cr
 \hfil#1#2\hfil\cr\noalign{\nointerlineskip}#1{\leftarrow}\mkern-6mu
 \cleaders\hbox{$#1\mkern-2mu{-}\mkern-2mu$}\hfill
 \mkern-6mu{-}\cr}}}
\let\amp@rs@nd@\relax
\newdimen\ex@
\ex@.2326ex
\newdimen\bigaw@
\newdimen\minaw@
\minaw@16.08739\ex@
\newdimen\minCDaw@
\minCDaw@2.5pc
\newif\ifCD@
\def\minCDarrowwidth#1{\minCDaw@#1}
\newenvironment{CD}{\@CD}{\@endCD}
\def\@CD{\def\A##1A##2A{\llap{$\vcenter{\hbox
 {$\scriptstyle##1$}}$}\Big\uparrow\rlap{$\vcenter{\hbox{%
$\scriptstyle##2$}}$}&&}%
\def\V##1V##2V{\llap{$\vcenter{\hbox
 {$\scriptstyle##1$}}$}\Big\downarrow\rlap{$\vcenter{\hbox{%
$\scriptstyle##2$}}$}&&}%
\def\={&\hskip.5em\mathrel
 {\vbox{\hrule width\minCDaw@\vskip3\ex@\hrule width
 \minCDaw@}}\hskip.5em&}%
\def\verteq{\Big\Vert&&}%
\def\noarr{&&}%
\def\vspace##1{\noalign{\vskip##1\relax}}\relax\let\amp@rs@nd@&\iffalse}\fi
 \CD@true\vcenter\bgroup\relax\let\\=\cr\iffalse}\fi\tabskip\z@skip\baselineskip20\ex@
 \lineskip3\ex@\lineskiplimit3\ex@\halign\bgroup
 &\hfill$\m@th##$\hfill\cr}
\def\@endCD{\cr\egroup\egroup}
\def\>#1>#2>{\amp@rs@nd@\setbox\z@\hbox{$\scriptstyle
 \;{#1}\;\;$}\setbox\@ne\hbox{$\scriptstyle\;{#2}\;\;$}\setbox\tw@
 \hbox{$#2$}\ifCD@
 \global\bigaw@\minCDaw@\else\global\bigaw@\minaw@\fi
 \ifdim\wd\z@>\bigaw@\global\bigaw@\wd\z@\fi
 \ifdim\wd\@ne>\bigaw@\global\bigaw@\wd\@ne\fi
 \ifCD@\hskip.5em\fi
 \ifdim\wd\tw@>\z@
 \mathrel{\mathop{\hbox to\bigaw@{\rightarrowfill}}\limits^{#1}_{#2}}\else
 \mathrel{\mathop{\hbox to\bigaw@{\rightarrowfill}}\limits^{#1}}\fi
 \ifCD@\hskip.5em\fi\amp@rs@nd@}
\def\<#1<#2<{\amp@rs@nd@\setbox\z@\hbox{$\scriptstyle
 \;\;{#1}\;$}\setbox\@ne\hbox{$\scriptstyle\;\;{#2}\;$}\setbox\tw@
 \hbox{$#2$}\ifCD@
 \global\bigaw@\minCDaw@\else\global\bigaw@\minaw@\fi
 \ifdim\wd\z@>\bigaw@\global\bigaw@\wd\z@\fi
 \ifdim\wd\@ne>\bigaw@\global\bigaw@\wd\@ne\fi
 \ifCD@\hskip.5em\fi
 \ifdim\wd\tw@>\z@
 \mathrel{\mathop{\hbox to\bigaw@{\leftarrowfill}}\limits^{#1}_{#2}}\else
 \mathrel{\mathop{\hbox to\bigaw@{\leftarrowfill}}\limits^{#1}}\fi
 \ifCD@\hskip.5em\fi\amp@rs@nd@}
\newenvironment{CDS}{\@CDS}{\@endCDS}
\def\@CDS{\def\A##1A##2A{\llap{$\vcenter{\hbox
 {$\scriptstyle##1$}}$}\Big\uparrow\rlap{$\vcenter{\hbox{%
$\scriptstyle##2$}}$}&}%
\def\V##1V##2V{\llap{$\vcenter{\hbox
 {$\scriptstyle##1$}}$}\Big\downarrow\rlap{$\vcenter{\hbox{%
$\scriptstyle##2$}}$}&}%
\def\={&\hskip.5em\mathrel
 {\vbox{\hrule width\minCDaw@\vskip3\ex@\hrule width
 \minCDaw@}}\hskip.5em&}
\def\verteq{\Big\Vert&}
\def\novarr{&}
\def\noharr{&&}
\def\SE##1E##2E{\slantedarrow(0,18)(4,-3){##1}{##2}&}
\def\SW##1W##2W{\slantedarrow(24,18)(-4,-3){##1}{##2}&}
\def\NE##1E##2E{\slantedarrow(0,0)(4,3){##1}{##2}&}
\def\NW##1W##2W{\slantedarrow(24,0)(-4,3){##1}{##2}&}
\def\slantedarrow(##1)(##2)##3##4{%
\thinlines\unitlength1pt\lower 6.5pt\hbox{\begin{picture}(24,18)%
\put(##1){\vector(##2){24}}%
\put(0,8){$\scriptstyle##3$}%
\put(20,8){$\scriptstyle##4$}%
\end{picture}}}
\def\vspace##1{\noalign{\vskip##1\relax}}\relax\let\amp@rs@nd@&\iffalse}\fi
 \CD@true\vcenter\bgroup\relax\let\\=\cr\iffalse}\fi\tabskip\z@skip\baselineskip20\ex@
 \lineskip3\ex@\lineskiplimit3\ex@\halign\bgroup
 &\hfill$\m@th##$\hfill\cr}
\def\@endCDS{\cr\egroup\egroup}
\newdimen\TriCDarrw@
\newif\ifTriV@
\newenvironment{TriCDV}{\@TriCDV}{\@endTriCD}
\newenvironment{TriCDA}{\@TriCDA}{\@endTriCD}
\def\@TriCDV{\TriV@true\def\TriCDpos@{6}\@TriCD}
\def\@TriCDA{\TriV@false\def\TriCDpos@{10}\@TriCD}
\def\@TriCD#1#2#3#4#5#6{%
\setbox0\hbox{$\ifTriV@#6\else#1\fi$}
\TriCDarrw@=\wd0 \advance\TriCDarrw@ 24pt
\advance\TriCDarrw@ -1em
\def\SE##1E##2E{\slantedarrow(0,18)(2,-3){##1}{##2}&}
\def\SW##1W##2W{\slantedarrow(12,18)(-2,-3){##1}{##2}&}
\def\NE##1E##2E{\slantedarrow(0,0)(2,3){##1}{##2}&}
\def\NW##1W##2W{\slantedarrow(12,0)(-2,3){##1}{##2}&}
\def\slantedarrow(##1)(##2)##3##4{\thinlines\unitlength1pt
\lower 6.5pt\hbox{\begin{picture}(12,18)%
\put(##1){\vector(##2){12}}%
\put(-4,\TriCDpos@){$\scriptstyle##3$}%
\put(12,\TriCDpos@){$\scriptstyle##4$}%
\end{picture}}}
\def\={\mathrel {\vbox{\hrule
   width\TriCDarrw@\vskip3\ex@\hrule width
   \TriCDarrw@}}}
\def\>##1>>{\setbox\z@\hbox{$\scriptstyle
 \;{##1}\;\;$}\global\bigaw@\TriCDarrw@
 \ifdim\wd\z@>\bigaw@\global\bigaw@\wd\z@\fi
 \hskip.5em
 \mathrel{\mathop{\hbox to \TriCDarrw@
{\rightarrowfill}}\limits^{##1}}
 \hskip.5em}
\def\<##1<<{\setbox\z@\hbox{$\scriptstyle
 \;{##1}\;\;$}\global\bigaw@\TriCDarrw@
 \ifdim\wd\z@>\bigaw@\global\bigaw@\wd\z@\fi
 \mathrel{\mathop{\hbox to\bigaw@{\leftarrowfill}}\limits^{##1}}
 }
 \CD@true\vcenter\bgroup\relax\let\\=\cr\iffalse}\fi
 \tabskip\z@skip\baselineskip20\ex@
 \lineskip3\ex@\lineskiplimit3\ex@
 \ifTriV@
 \halign\bgroup
 &\hfill$\m@th##$\hfill\cr
#1&\multispan3\hfill$#2$\hfill&#3\\
&#4&#5\\
&&#6\cr\egroup%
\else
 \halign\bgroup
 &\hfill$\m@th##$\hfill\cr
&&#1\\%
&#2&#3\\
#4&\multispan3\hfill$#5$\hfill&#6\cr\egroup
\fi}
\def\@endTriCD{\egroup}
\newcommand{\sA}{{\mathcal A}}
\newcommand{\sB}{{\mathcal B}}
\newcommand{\sC}{{\mathcal C}}
\newcommand{\sD}{{\mathcal D}}
\newcommand{\sE}{{\mathcal E}}
\newcommand{\sF}{{\mathcal F}}
\newcommand{\sG}{{\mathcal G}}
\newcommand{\sH}{{\mathcal H}}
\newcommand{\sI}{{\mathcal I}}
\newcommand{\sJ}{{\mathcal J}}
\newcommand{\sK}{{\mathcal K}}
\newcommand{\sL}{{\mathcal L}}
\newcommand{\sM}{{\mathcal M}}
\newcommand{\sN}{{\mathcal N}}
\newcommand{\sO}{{\mathcal O}}
\newcommand{\sP}{{\mathcal P}}
\newcommand{\sQ}{{\mathcal Q}}
\newcommand{\sR}{{\mathcal R}}
\newcommand{\sS}{{\mathcal S}}
\newcommand{\sT}{{\mathcal T}}
\newcommand{\sU}{{\mathcal U}}
\newcommand{\sV}{{\mathcal V}}
\newcommand{\sW}{{\mathcal W}}
\newcommand{\sX}{{\mathcal X}}
\newcommand{\sY}{{\mathcal Y}}
\newcommand{\sZ}{{\mathcal Z}}
\newcommand{\A}{{\mathbb A}}
\newcommand{\B}{{\mathbb B}}
\newcommand{\C}{{\mathbb C}}
\newcommand{\D}{{\mathbb D}}
\newcommand{\E}{{\mathbb E}}
\newcommand{\F}{{\mathbb F}}
\newcommand{\G}{{\mathbb G}}
\newcommand{\BH}{{\mathbb H}}
\newcommand{\I}{{\mathbb I}}
\newcommand{\J}{{\mathbb J}}
\renewcommand{\L}{{\mathbb L}}
\newcommand{\M}{{\mathbb M}}
\newcommand{\N}{{\mathbb N}}
\newcommand{\BP}{{\mathbb P}}
\newcommand{\Q}{{\mathbb Q}}
\newcommand{\R}{{\mathbb R}}
\newcommand{\BS}{{\mathbb S}}
\newcommand{\T}{{\mathbb T}}
\newcommand{\U}{{\mathbb U}}
\newcommand{\V}{{\mathbb V}}
\newcommand{\W}{{\mathbb W}}
\newcommand{\X}{{\mathbb X}}
\newcommand{\Y}{{\mathbb Y}}
\newcommand{\Z}{{\mathbb Z}}
\newcommand{\id}{{\rm id}}
\newcommand{\rk}{{\rm rank}}
\newcommand{\pr}{{\rm pr}}
\newcommand{\roundup}[1]{\ulcorner{#1}\urcorner}
\newcommand{\rounddown}[1]{\llcorner{#1}\lrcorner}
\newcommand{\END}{{\mathbb E}{\rm nd}}
\newcommand{\End}{{\rm End}}
\newcommand{\Hg}{{\rm Hg}}
\newcommand{\tr}{{\rm tr}}
\newcommand{\Sl}{{\rm Sl}}
\newcommand{\Gl}{{\rm Gl}}
\newcommand{\Sp}{{\rm Sp}}
\newcommand{\MT}{{\rm MT}}
\newcommand{\Cor}{{\rm Cor}}
\newcommand{\Hom}{{\sH}{\rm om}}
\newcommand{\Mon}{{\rm Mon}}
\newcommand{\s}{{\rm sl}}
\newcommand{\ch}{{\rm c}}
\newcommand{\gr}{{\mathfrak g \mathfrak r}}
\newcommand{\Gr}{{\G{\rm r}}}
\begin{abstract}
Given an open subset $U$ of a projective curve $Y$ and a
smooth family $f:V\to U$ of curves, with semi-stable reduction over $Y$,
we show that for a sub variation $\V$ of Hodge structures
of $R^1f_*\C_V$ with $\rk(\V)>2$ the Arakelov inequality must be strict. 
For families of $n$-folds we prove a similar result under the assumption that
the $(n,0)$ component of the Higgs bundle of $\V$ defines a birational map.
\end{abstract}
\maketitle
\section*{Introduction}
\renewcommand{\thethm}{\arabic{thm}} \renewcommand{\theequation}{\arabic{equation}} 
Let $f:X\to Y$ be a complex-projective family of $n$-folds over a curve $Y$, i.e. $X$ is a smooth 
complex-projective $n+1$-fold, $Y$ is a complex projective curve, and $f$ is a surjective
morphism with connected fibres. Let $S\subset Y$ be the critical locus of $f$, so that $f$ restricts to a 
smooth projective morphism 
$$
f:V=f^{-1}(U)\>>> U = Y\setminus S.
$$
The direct image sheaf $R^nf_*\C_V$ is a local system on $U$ which underlies a variation of Hodge structures $\V$ of weight $n$. 
By the monodromy theorem the local monodromy operators around the points of $S$ are quasi-unipotent,
so replacing $Y$ by a suitable ramified covering one may assume that these are all unipotent.
In this note we assume however that this is already the case for $f:X\to Y$,
or the slightly stronger condition that $f$ is semi-stable, hence all fibres $f^{-1}(y)$
reduced normal crossing divisors. 

As a motivation of the paper consider the case $n=1$, i.e. a family of semi-stable
curves of genus $g$. It is shown in \cite{Bea1} that for a non-isotrivial family of curves over 
$Y=\BP^1$ one has $\#S \geq 4$ and that $\#S =4$ implies that the irregularity of the total 
space $X$ is zero. By \cite{Tan} (see also \cite{Ngu}), for $g\geq 2$,
one has $\#S \geq 5$. This remarkable result follows from
Beauville's observation and from Tan's strict Arakelov inequality  
\begin{equation}\label{eqin.1}
\deg(f_*\omega_{X/Y}) < \frac{1}{2}\cdot g \cdot\deg(\Omega_Y^1(\log S)).
\end{equation}
Although the bound $\#S\geq 5$ is optimal, we will show in this note, that
the inequality (\ref{eqin.1}) can be strengthened, and that under certain assumptions it  
extends to the case $n>1$.

Since $Y$ is a curve, the vector bundle $\V\otimes \sO_U$ extends to a bundle $\sH$ on $C$ in such a way that the Gauss-Manin connection acquires logarithmic singularities. We choose for $\sH$ the Deligne extension, i.e. an extension such that the real part of the local residues are zero.
The Hodge filtration extends to a holomorphic filtration on $\sH$. The extended Gauss-Manin connection
defines on the associated graded bundle the structure of a so-called Higgs bundle $(F,\tau)$,
i.e. a collection of vector bundle maps
$$
\tau_{p,q}:F^{p,q}\>>> F^{p-1,q+1}\otimes \Omega_Y^1(\log S), \ \ \ \ \ \ p+q=n.
$$
The component $F^{0,n}$ can be identified with $f_*\omega_{X/Y}$. In particular $F^{0,n}$ is 
numerically effective. On the Higgs-bundle side we can use \cite{VZ3} which bounds the degree
of $F^{0,n}$ (see Lemma \ref{in.4}). 

Because of Koll\'ar's decomposition \cite{Kol}
$$
f_*\omega_{X/Y}=\sA\oplus\sB,
$$
into an ample sheaf $\sA$ and a flat subsheaf $\sB$, we have $\deg(F^{0,n})=\deg(\sA)$,
and the main result can be stated using this bundle:

\begin{theorem}\label{in.1}
Keeping the notations introduced above, let $f:X \to Y$ be a semi-stable
non-isotrivial family of $n$-folds. Assume that either
\begin{enumerate} 
\item[a.] $f^*\sA \to \omega_{X/Y}$ defines a birational $Y$-morphism
$\eta:X \to \BP(\sA)$,
\item[b.] or that $n=1$ and $\rk(\sA)\geq 2$.
\end{enumerate}
Then
$$
\deg(\sA)=\deg(f_*\omega_{X/Y}) < \frac{n}{2}\cdot\rk(\sA)\cdot\deg(\Omega_Y^1(\log S)).
$$
\end{theorem}

Let us consider for a moment the corresponding question for semi-stable families 
$h:A\to Y$ of $g$-dimensional Abelian varieties, smooth over $U$. 
Here one considers the weight one variation of Hodge structures $R^1h_*\C_{g^{-1}(U)}$. So the $(1,0)$
part of the Higgs bundle is $F^{1,0}=h_*\Omega^1_{A/Y}(\log h^{-1}S)$. Again one has Koll\'ar's decomposition $F^{1,0}=\sA\oplus\sB$ into an ample sheaf $\sA$ and a flat sheaf $\sB$. As shown in \cite{VZ2}, one has again an in equality
\begin{equation}\label{eqin.2}
\deg(\sA)\leq \frac{1}{2}\cdot\rk(\sA)\cdot\deg(\Omega_Y^1(\log S)),
\end{equation}
and the equality in (\ref{eqin.2}) 
implies that $U$ is a Shimura curve of Hodge type, as defined in \cite{Mum}. By definition $U$ is an \'etale covering of the moduli space of Abelian varieties with prescribed Mumford-Tate group, and 
$V\to U$ is the universal family. The classification of Shimura curves, or the explicite description of the
variation of Hodge structures $R^1h_*\C_{g^{-1}(U)}$ in  \cite{VZ2} imply that
$\sA$ is a direct sum of $\rk(\sA)$ copies of an invertible sheaf if $S\neq \emptyset$, and
that up to constant factors $h:A\to Y$ is isogenous to the selfproduct of $\rk(\sA)$ copies of a modular family of elliptic curves. On the other hand, for $S=\emptyset$ the sheaf $\sA$ can not contain an invertible direct factor if (\ref{eqin.2}) is an equality.

If $f:X\to Y$ is a semi-stable family of curves, smooth over $U$, and with $g:A\to Y$ as family of Jacobians, one has 
$$
f_*\omega_{X/Y}=h_*\Omega^1_{A/Y}(\log h^{-1}S)=\sA\oplus\sB.
$$ 
As we recall below the characterization of Teichm\"uller curves in \cite{Moe} implies that
equality in (\ref{eqin.2}) can not hold if $\sA$ is a direct sum of two or more
invertible sheaves. In different terms, if the moduli space $M_g$ of non-singular 
curves of genus $g$ contains a non-compact Shimura curve, the moving part $\sA$ of the
Higgs bundle has to be of rank one. As shown in \cite{Moe2} such ``Shimura-Teichm\"uller curves''
in $M_g$ can only exist for $g=3$. 

To prove Theorem \ref{in.1}, b), it remains to exclude the case
$S=\emptyset$, or equivalently the case where $\sA$ does not contain
an invertible direct factor. Using the characterization of Shimura curves one obtains:

\begin{corollary}[\cite{Moe2} for $S\neq \emptyset$]\label{in.2}
The moduli space $M_g$ of non-singular curves of genus $g$ does not contain a 
compact Shimura curve, and for $g\neq 3$ it does not contain any
Shimura curve at all.
\end{corollary}
\cite{Moe2} gives an explicite example of a non-compact Shimura curve in $M_3$.
So in Theorem \ref{in.1} it is not sufficient to assume that $\sA\neq 0$, even in case the
general fibre of $f$ is of general type. We hope that the condition ``$\eta$ birational'' in Theorem \ref{in.1}, a), can be replaced by ``$\eta$ generically finite'', 
but we were unable to prove Theorem \ref{in.1}, a), under this assumption.

Let us emphasize that the Corollary \ref{in.2} is a nice answer to the wrong question.
For $g$ sufficiently large, there should not exist Shimura curves in the closure of $M_g$ in $A_g$,
but the methods presented here do not allow any result in this direction.\\

Theorem \ref{in.1} generalizes to other decompositions of $f_*\omega_{X/Y}$.
To explain the set-up, recall that Koll\'ars decomposition is induced
by a decomposition of variations of Hodge structures
$$
R^nf_*\C_{V}=\V\oplus \U
$$
with $\U$ unitary and where the $(n,0)$-part of the logarithmic Higgs bundle of $\V$ is $\sA$. 
Dropping the first condition one can consider any sub variation of Hodge structures $\V$ of
$R^nf_*\C_{V}$. For $n>1$ we will pose the condition that the $(n,0)$ part $F^{n,0}$ of the Higgs bundle of $\V$ is non-isotrivial in the sense that the image of the rational map
$$
\varphi_\V:X\>>> \BP(F^{n,0})
$$
induced by the natural evaluation map $f^*F^{n,0} \to \omega_{X/Y}$ does not produce
an isotrivial family over $Y$. This hypothesis is verified if $f$ is non-isotrivial and
$\varphi_\V$ itself birational.

Abusing notations let 
$$
(F,\tau)=\big(\bigoplus_{p+q=n}F^{p,q},\bigoplus_{p+q=n}\tau_{p,q}\big)
$$
be the logarithmic Higgs bundle associated to the sub variation $\V$. The bundle maps $\tau_{p,q}$
can be iterated to obtain maps
$$
\tau^{(\ell)}: F^{n,0} \>>> F^{n-1,1}\otimes \Omega^1_Y(\log S)\>>> \cdots
\>>> F^{n-\ell,\ell}\otimes S^\ell(\Omega^1_Y(\log S)).
$$
\begin{definition}\label{in.3} \ 
\begin{enumerate}
\item[i.] We call $\tau^{(n)}:F^{n,0} \to F^{0,n}\otimes S^n(\Omega^1_Y(\log S))$
the Griffiths-Yukawa coupling of $\V$ (or of $f$ in case $\V=R^nf_*\C_V$).
\item[ii.] The Griffiths-Yukawa coupling of $\V$ is maximal, if 
$F^{n,0}\neq 0$ and if $\tau^{(n)}$ is an isomorphism.
\item[iv.] The Higgs field is strictly maximal, if $F^{0,n}\neq 0$ and if all the $\tau_{p,q}$
are all isomorphisms.
\end{enumerate}
\end{definition}
In Section \ref{ar} we will show:
\begin{lemma}\label{in.4}
Assume that $\V$ is a non-trivial variation of polarized complex Hodge structures 
of weight $n$ with unipotent local monodromy in all $s\in S$, and with
logarithmic Higgs bundle $(\bigoplus F^{p,q},\tau_{p,q})$. Then 
\begin{equation}\label{eqin.3}
\deg(F^{n,0}) \leq \frac{n}{2}\cdot\rk(F^{n,0})\cdot\deg(\Omega^1_Y(\log S)),
\end{equation}
and (\ref{eqin.3}) is an equality if and only if the Griffiths-Yukawa coupling
of $\V$ is maximal. Moreover, in this case one has a decomposition $\V=\V_1\oplus \V_2$ where
the Higgs field of $\V_1$ is strictly maximal and where
$\V_2$ is a variation of polarized complex Hodge structures, zero in bidegree
$(n,0)$. 
\end{lemma}
If $\#S$ is even, \cite[3.4]{VZ3} gives a more precise description of
$\V_1$. Choose a logarithmic theta characteristic, 
i.e. an invertible sheaf $\sL$ with $\sL^2=\Omega_Y^1(\log S)$, and write 
$\L$ for the local system with Higgs bundle $\sL\oplus \sL^{-1}$ and Higgs field
$$
\sL \>{\rm id}>> \sL^{-1}\otimes\Omega_Y^1(\log S). 
$$
$\L$ is unique up to the tensor product with a unitary rank one local system,
induced by a two-division point of ${\rm Pic}^0(Y)$.
\begin{addendum}\label{in.5}
If $\#S$ is even, and $\L$ induced by a logarithmic theta characteristic, then there exists 
a unitary local system $\T$, regarded as a variation of Hodge structures in bidegree $(0,0)$,
with $\V_1=S^n(\L)\otimes \T$. 
\end{addendum}
The assumption that the local monodromy operators are unipotent or that the family $f:X\to Y$
is semistable is not really needed at this point. Without it the $(n,0)$ part of 
$R^nf_*\C_V$ is $f_*\Omega^n_{X/Y}(\log f^{-1}(S))$. This sheaf can only become larger under semistable reduction. For a subsheaf $F^{n,0}$ (\ref{eqin.3}) can only be an equality,
if $F^{n,0}$ is compatible with passing to a semi-stable model over some covering of $Y$, \'etale
over $U$. Hence one does not lose any information working with such a model from the start.

In view of Lemma \ref{in.4} Theorem \ref{in.1} is a special case of the following Theorem: 
\begin{theorem}\label{in.7} 
Let $f:X\to Y$ be a family of $n$-folds over a curve $Y$, semi-stable or with unipotent monodromy
in $s\in S$, and let $V \to U$ be the smooth part of $f$.
Let $\V$ be a complex polarized sub variation of Hodge structures of $R^nf_*\C_V$ with 
logarithmic Higgs bundle $(F,\tau)$. Assume one of the following conditions:
\begin{enumerate}
\item[i.] The $(n,0)$ component $F^{n,0}$ of $\V$ is non-isotrivial.
\item[ii.] $f$ is a family of curves, and $\rk(F^{1,0})>1$. 
\end{enumerate}
Then the Griffiths-Yukawa coupling of $\V$ is not maximal.
\end{theorem}
In \cite{VZ3} we studied families of higher dimensional minimal manifolds with a  maximal
Griffiths-Yukawa coupling for $R^nf_*\C_V$. We asked, whether such families can exist
if the Kodaira dimension of the general fibre $F$ is positive, or if $p_g(F)>1$. Theorem \ref{in.7}, i),
gives a negative answer, under the stronger assumption that the canonical linear system
$|\omega_F|$ defines a birational map. Again, we hope that Theorem \ref{in.7}, i), also holds true
if $|\omega_F|$ defines a generically finite map. Part ii) confirms the latter  
for families of curves. Here two additional facts are needed. The first one, proved in Section \ref{mu}
states that for $n=1$ the isotriviality of $F^{1,0}$ and the strict maximality of the Higgs field 
imply that $F^{1,0}$ is a direct sum of line bundles. 
The second one comes from the theory of Teichm\"uller curves (see \cite{McM1}). 

Roughly speaking, one considers geodesics in the Teichm\"uller space,
constructed by an $\Sl(2,\R)$-action on the real and imaginary part of a given holomorphic differential form. If the quotient by a suitable lattice in $\Sl(2,\R)$ is an algebraic curve, it is called a 
Teichm\"uller curve. We will not need this definition, since M. M\"oller has given in \cite{Moe}
an algebraic characterization of such curves:
\begin{theorem}[\cite{Moe}]\label{in.8}
A semi-stable family of curves $f:X\to Y$ of genus $g$ defines a Teichm\"uller curve in the moduli space $M_g$ if and only if $R^1f_*\C_V$ contains a sub variation of Hodge structures $\V$ of rank two
with a strictly maximal Higgs field.
\end{theorem}
Using the description of such curves in \cite{McM1} one finds that there is no other sub variation $\V$ 
with a strictly maximal Higgs field, and that $S\neq \emptyset$. This allows in the proof of Theorem
\ref{in.7}, ii), to assume that $F^{1,0}$ is not the direct sum of line bundles, or equivalently that
$S\neq\emptyset$. In Section \ref{mu} we will see that this assumption, as well as the assumption on the
non-isotriviality in part i), allows to use the maximality of the Griffiths-Yukawa coupling 
to construct a section (replacing $U$ by an \'etale covering) of the sheaf
$\omega_{X/Y}^\nu\otimes f^*\omega_Y(S)^{-\frac{n\cdot\nu}{2}}$, vanishing with order $\nu$ along some fibre. In the last section we will show, that such a section can not exist.\\

Let us return to semi-stable families over $\BP^1$. It is easy to construct a family of elliptic curves $E\to \BP^1$ with three multiple fibres, two of which are semi-stable, and the third one 
has semi-stable reduction over a covering of degree two. In fact, one just has to take a two-fold covering of 
$\BP^1\times \BP^1$, ramified over 
$$
\{0\} \times \BP^1 +\{1\}\times \BP^1+ \{\infty\} \times \BP^1+\Delta + \BP^1 \times \{0\}.
$$
Then the monodromy of $Z=E\times_{\BP^1} E \to \BP^1$ is unipotent,
Dividing by the involution one obtains a family of $K3$ surfaces $f:X\to \BP^1$ with $3$ singular fibres
and with unipotent local monodromy operators. The family is non-isotrivial, hence the $(2,0)$ component of its Higgs bundle has positive degree. The inequality (\ref{eqin.3}) implies that it is $\sO_{\BP^1}(1)$,
and the Griffiths-Yukawa-coupling has to be maximal (see \cite[Example 7.5]{VZ3} for similar calculations).  
With a little bit of work, one can presumably show that this family has a birational model 
with $3$ singular semi-stable fibres. In odd dimensions similar examples can not exist: 
  
\begin{proposition}\label{in.9}
Let $f:X\to \BP^1$ be a non isotrivial family of $n$-dimensional varieties
with general fibre $F$. Assume that the local monodromy operators in $s\in S$ are uni-potent and that $n$ is odd. If $R^nf_*\C_V$ contains a non trivial local sub system with a maximal Griffiths-Yukawa
coupling, $\# S$ is even, hence $\#S \geq 4$. 
\end{proposition}

For $n\geq 3$ and odd, Proposition \ref{in.9} does not allow to improve the known bound $\#S \geq 3$
for the number of singular fibres. Contrary to the case of curves (or surfaces), the existence of a family with  $4$ (or $3$) singular fibres does not imply the maximality of the Griffiths-Yukawa coupling,
even if one assumes that the local Torelli theorem holds.\\
\ \\
Martin M\"oller introduced us to the theory of Teichm\"uller curves. 
His article \cite{Moe}, and his examples of special Teichm\"uller curves in \cite{Moe2} were of 
high interest for our understanding of sub systems with a maximal Griffiths-Yukawa coupling, and 
part ii) of Theorem \ref{in.7} is a consequence of \cite{Moe}.
Although this aspect does not appear in the article, its starting point was
a try to understand the relation between geodesic curves in moduli spaces for different 
natural metrics in $M_g$ and $A_g$, in particular a long discussion with 
Shing-Tung Yau. We would like to thank both of them for their 
interest and help. We also thank the referee of a first version of the article,
for pointing out several ways to improve its presentation.

This note grew out of discussions started when the first named author visited
the Institute of Mathematical Science and the Department of Mathematics at the Chinese
University of Hong Kong. The final steps were done when he visited the IAS, Princeton.
He would like to thank the members of both Institutes for their hospitality.
\section{Arakelov inequalities}
\label{ar}
\renewcommand{\thethm}{\arabic{section}.\arabic{thm}} \renewcommand{\theequation}{\arabic{section}.\arabic{equation}}
Let us recall the Arakelov inequalities shown in Section 2 of \cite{VZ3}. Let $Y$ be a 
projective curve, $S$ a finite set of points, and let
$$
(F,\tau)=\big(\bigoplus_{p+q=n} F^{p,q},\bigoplus_{p+q=n} \tau_{p,q}\big)
$$
be the Higgs field corresponding to the Deligne extension of a polarized complex 
variation of Hodge structures $\V$ of weight $n$ on $Y\setminus S$
with unipotent local monodromy operators around the points in $S$. We write
$$ F_0^{p,q}={\rm ker}(\tau_{p,q}: F^{p,q}\to F^{p-1,q+1}\otimes
\Omega^1_Y(\log S)),\quad h^{p,q}_0={\rm rk} F_0^{p,q}.$$
\begin{proposition}\label{ar.1} \
\begin{enumerate}
\item[a.] One has
$$
\deg(F^{n,0}) \leq \frac{n}{2}\cdot(h^{n,0}-h_0^{n,0})
\cdot \deg(\Omega^1_Y(\log S)).
$$
\item[b.] Assume that $\tau_{n,0}\neq 0$. Then $0 < \deg(F^{n,0})$
and $\deg(\Omega^1_Y (\log S))>0$. In particular, for $Y=\BP^1$ one has $\#S\geq 3$. 
\item[c.] Let $E^{n,0}$ be a subsheaf of $F^{n,0}$ such that 
$$
\deg(E^{n,0}) = \frac{n}{2}\cdot\rk(E^{n,0})\cdot \deg(\Omega_Y^1(\log S)).
$$
Then one has a decomposition of polarized variations of Hodge structures
$\V=\W\oplus\W'$, and subsheaves $E^{p,n-p}$ of $F^{p,n-p}$, for $p=n-1,\ldots,0$,
such that the Higgs bundle of $\W$ is 
$$
(E,\theta)=\big(\bigoplus_{p+q=n}E^{p,q},\bigoplus_{p+q=n}\theta_{p,q}\big),
$$ 
and the Higgs field $\theta$ is strictly maximal.
\end{enumerate}
\end{proposition}
\begin{proof}
a) and b) are special cases of \cite[Proposition 2.1]{VZ3}. For c)
consider the Higgs sub bundles $E=\bigoplus_{p+q=n}E^{p,q}$ of $F$ with
$$
E^{n-\ell,\ell}=\tau^{(\ell)}(E^{n,0})\otimes \Omega^1_Y(\log S)^{-\ell}.
$$
Assume the equation in c) holds. Let $K^{n-\ell,\ell}$ be the kernel of
$$
\tau_{n-\ell,\ell}|_{E^{n-\ell,\ell}}:E^{n-\ell,\ell} \>>> 
E^{n-\ell-1,\ell+1}\otimes\Omega^1_Y(\log S)\subset
F^{n-\ell-1,\ell+1}\otimes\Omega^1_Y(\log S).
$$
Since $K^{n-\ell,\ell}$ is a sub Higgs bundle, by Simpson \cite{Sim} 
$\deg(K^{n-\ell,\ell})\leq 0$. One obtains 
\begin{multline*}
\deg(E^{n-\ell-1,\ell+1}) = \deg(\tau_{n-\ell,\ell}(E^{n-\ell,\ell})) - 
\rk(E^{n-\ell-1,\ell+1})\cdot
\deg(\Omega_Y^1(\log S))\\
\geq \deg(E^{n-\ell,\ell}) - \rk(E^{n-\ell-1,\ell+1})\cdot
\deg(\Omega_Y^1(\log S)).
\end{multline*}
Then
$$
\deg(E^{n-\ell-1,\ell+1}) \geq \deg(E^{n,0}) - \sum_{k=1}^{\ell+1}\rk(E^{n-k,k})\cdot
\deg(\Omega_Y^1(\log S)),
$$
and adding up one finds by assumption
\begin{multline*}
\deg(E) \geq (n+1)\cdot \deg(E^{n,0}) - \sum_{\ell=0}^{n} 
\sum_{k=1}^{\ell+1}\rk(E^{n-k,k})\cdot \deg(\Omega_Y^1(\log S))\\
=\big((n+1)\cdot\frac{n}{2}\cdot \rk(E^{n,0}) - \sum_{k=1}^n (n-k+1)\cdot \rk(E^{n-k,k})
\big)\cdot \deg(\Omega_Y^1(\log S))\geq 0.
\end{multline*}
$E$ is a sub Higgs bundle of $F$, hence $\deg(E)\leq 0$, and all those inequalities are
equalities. For the last one this implies that $\rk(E^{n,0})=\rk(E^{n-k,k})$, for all $k$
and 
$$
\theta_{p,q}=\tau_{p,q}|_{E^{p,q}}:E_{p,q} \>>> E_{p-1,q+1}\otimes\Omega_Y^1(\log S)=
\tau_{p,q}(E^{p,q})
$$
must be an isomorphism. Moreover $\deg(E)=0$ implies that $E$ is a direct factor of 
$(F,\tau)$, and that it is the Higgs field of a sub local system $\W$ of $\V$.
Then $\W$ is a polarizable $\C$-variation of Hodge structures and by \cite{Del}
one obtains $\V=\W\oplus \W'$.
\end{proof}
\begin{proof}[Proof of Lemma \ref{in.4} and Addendum \ref{in.5}]
Applying \ref{ar.1}, c), to $E^{n,0}=F^{n,0}$ one sees that
equality in (\ref{eqin.3}) implies the maximality of the Griffiths-Yukawa coupling
and the existence of the decomposition $\V=\V_1\oplus\V_2$.

On the other hand, if the Griffiths-Yukawa coupling is maximal one has
$$
(F^{n,0})^\vee=F^{0,n}\cong F^{n,0}\otimes \Omega_Y^1(\log S)^{-n}
$$
and $-\deg(F^{n,0})=\deg(F^{n,0}) - \rk(F^{n,0})\cdot n\cdot \deg(\Omega^1_Y(\log S))$,
as claimed. If $\#S$ is even the tensor product decomposition $\V_1=S^n(\L)\otimes \T$ 
has been obtained in \cite[3.4]{VZ3}. 
\end{proof}
\begin{proof}[Proof of Proposition \ref{in.9}]
By assumption $R^nf_*\C_V$ contains a sub variation of Hodge structures $\V$ with a 
maximal Griffiths-Yukawa coupling, and \ref{ar.1}, c), allows to assume that the Higgs field
of $\V$ is strictly maximal. 
Let $(\bigoplus F^{p,n-p},\tau)$ be the Higgs bundle of $\V$. One can write
$$
F^{0,n}=\bigoplus_{i=1}^{r}\sO_{\BP^1}(\nu_i) \mbox{ \ \ with \ \ }
\nu_1\geq \nu_2 \geq \cdots \geq \nu_r.
$$
If $\#S$ is odd, consider a twofold covering $\varphi:Y'=\BP^1 \to \BP^1$ ramified in exactly 
two points in $S$. Write $S'=\tau^{-1}(S)$ and $f':X'\to Y'$ for the pullback family.
The pullback $\V'$ of $\V$ is a sub variation of Hodge structures of $R^nf'_*\C_{V'}$, with 
Higgs field
$$
(\bigoplus_{p=0}^n {F'}^{p,n-p}=\bigoplus_{p=0}^n \varphi^*{F}^{p,n-p} ,\tau').
$$
Obviously the Higgs field $\tau'$ is still maximal, and the Addendum \ref{in.5} implies that
$\V_1=S^n(\L)\otimes \T$, for a rank two variation of Hodge structures $\L$ and for a unitary
bundle $\T$ in bidegree $(0,0)$. If $\#S$ is even, the same holds true on $\BP^1$ itself.

In both cases one obtains $\nu_1=\nu_2=\cdots=\nu_r=\nu$ and the Arakelov equality reads
\begin{gather*}
r\cdot \nu=\deg({F}^{n,0}) = \frac{n\cdot r}{2}\cdot \deg(\Omega^1_{\BP^1}(\log S)),\\
\mbox{ \ \ or \ \ }
\nu = \frac{n}{2}\cdot \deg(\Omega^1_{\BP^1}(\log S))=\frac{n}{2}\cdot (-2+\#S) .
\end{gather*}
So the right hand side must be an integer, and for $n$ odd $\#S$ must be even.
\end{proof}

\section{The multiplication map}\label{mu}

In order to prove Theorem \ref{in.7} we have to exclude the existence of certain 
families $f:X\to Y$ of $n$-folds whose variation of Hodge structures $R^nf_*\C_V$ contains a
sub variation of complex polarized Hodge structures $\V$, with unipotent local monodromy operators, 
and with a strictly maximal Higgs field. Again
$$
\big(F= \bigoplus_{q=0}^n F^{n-q,q}, \ \tau \big) 
$$
denotes the Higgs field of $\V$.
\begin{proposition}\label{mu.1}
Assume that $F^{n,0}$ is non-isotrivial. Then there exists a finite covering $\varphi:\hat{Y}\to Y$, \'etale over $U$, 
such that for the induced family $\hat{f}:\hat{X}\to \hat{Y}$ with $\hat{X}$ a desingularization of
$X\times_Y\hat{Y}$ the following condition holds true:
\begin{enumerate}
\item[($*$)] For some $\nu$ sufficiently large and divisible by $2$
and for a point $y\in \hat{Y}$ in general position there exists a non-trivial section
of 
$$
{\hat{f}}_*\omega_{\hat{X}/\hat{Y}}^\nu\otimes \omega_{\hat{Y}}(\varphi^{-1}(S))^{-\frac{n\cdot\nu}{2}} \otimes \sO_{\hat{Y}}(-\nu\cdot y).
$$
\end{enumerate}
\end{proposition}
\begin{proposition}\label{mu.2}
If $f:X\to Y$ is a family of curves and $\rk(F^{1,0})\geq 2$.
Then there exists a finite covering $\varphi:\hat{Y}\to Y$, \'etale over $U$, such that either the condition {\rm($*$)} in \ref{mu.1} holds true, or 
\begin{enumerate}
\item[($**$)]
$\varphi^*F^{1,0}=\bigoplus^\ell \hat{\sL}$ for a line bundle $\hat{\sL}$ with $\hat{\sL}^2=\omega_{\hat{Y}}(\varphi^{-1}(S))$. 
\end{enumerate}
\end{proposition}
For both Propositions we start with the same construction. Since $\bar{\V}$ also has a strictly maximal Higgs field, we may enlarge $\V$ and assume that it is invariant under complex conjugation.

By assumption $\rk(F^{n,0})=\ell \geq 2$, hence $\tau_{n,0}\neq 0$ and Proposition \ref{ar.1} implies that
$2\cdot g(Y) -2 +\#S >0$. Then there exist coverings $\varphi:\hat{Y}\to Y$, \'etale over $U$, of arbitrarily 
high degree. In particular, dropping the $\hat{ \ \ }$ we may replace $U$ by an \'etale
covering, assume that $\#S$ is even, and choose a logarithmic theta characteristic $\sL$, i.e.
an invertible sheaf with $\sL^{2}=\Omega_Y^1(\log S)$.

The Addendum \ref{in.5} implies that $\V$ is of the form $S^n(\L)\otimes \T$, and since $\L$ is invariant
under complex conjugation, the same holds true for $\T$. So the Higgs field of $\V$ is of the form 
$$
(F=\bigoplus_{p+q=n}{F}^{p,q},\tau),
$$
with $F^{n,0}=\sL^n\otimes \sT$ and $F^{0,n}=\sL^{-n}\otimes \sT$ for $\sT$ unitary of rank $\ell$.
In particular, $\sT$ is self dual, and $\det(\sT)^2=\sO_Y$. Replacing $\sL$ and $\sT$ by 
$$
\det(\sT)^{\frac{1}{n\cdot\ell}}\otimes \sL \mbox{ \ \ and \ \ }
\det(\sT)^{\frac{-1}{\ell}}\otimes \sT
$$ 
one may assume that $\det(\sT)=\sO_Y$. The evaluation map
$$
\rho:f^*(\sT) \> >> \omega_{X/Y}\otimes f^*\sL^{-n}
$$
induces a rational map $\varrho: X \to \BP(\sT)$ over $Y$. 
Choose a blowing up $\delta:Z \to X$ such that $\varrho\circ \delta$ is a morphism, and consider the diagram
$$
\begin{TriCDV}
{Z}{ \> \sigma >>}{\BP(\sT)}
{\SE h EE}{\SW W \pi W}
{Y.}
\end{TriCDV}
$$
We write $W$ for the image of $\sigma$, 
$$
\pi'=\pi|_W:W\>>> Y \mbox{ \ \ and \ \ } \sO_{W}(\nu)=\sO_{\BP}(\nu)|_W.
$$ 
Then
$$
\sigma^*\sO_{\BP}(1)\subset \delta^*\omega_{X/Y}\otimes h^*\sL^{-n},
$$
and the multiplication map is given by
$$
S^\nu(\sT) = \pi_*\sO_{\BP}(\nu) \>m_\nu >> \pi'_*\sO_{W}(\nu) \> \subset >>
(h_*\delta^*\omega_{X/Y}^\nu)\otimes \sL^{-n\cdot\nu}.
$$
$\sT_\nu\subset \pi'_*\sO_{W}(\nu)$ will denote the image of $m_\nu$. For $\nu$
sufficiently large one has $\sT_\nu = \pi'_*\sO_{W}(\nu)$. 

\begin{claim}\label{mu.3}
If $\deg(\sT_\nu)=0$, for all $\nu \gg 1$, the fibre 
$$
\pi'^{-1}(y)\subset \pi^{-1}(y)\cong\BP^{\ell-1}
$$
is independent of $y\in Y$ up to the choice of coordinates in $\BP^{\ell-1}$. 
\end{claim}
\begin{proof}
The locally free sheaf $\sT$, as well as its symmetric products are poly-stable
of degree $0$. Let $\sK_\nu$ be the kernel of $m_\nu$. If $\deg(\sT_\nu)=0$ 
one obtains a decomposition $S^\nu(\sT)=\sK_\nu\oplus \sT_\nu$ with $\sK_\nu$
poly-stable of degree $0$, hence flat and unitary. As in \cite[p. 396]{Kol} 
this implies the isotriviality of $W$. To see that all fibres of $\pi'$ are really isomorphic,
we argue in a similar way, along the line of the proof of \cite[Theorem 4.33 or 4.34]{Vie}. 

Consider the projective bundle $\eta:\Pi=\BP({\sT^\vee}^{\oplus \ell})\to Y$. As in
\cite[Section 4.4]{Vie} for an effective divisor $\Delta$ the open subscheme $\Pi_0=\Pi\setminus \Delta$ 
is a $\BP{\rm Gl}(\ell,\C)$ torsor over $Y$, and 
$\sO_{\Pi}(\Delta)=\sO_{\Pi}(\ell)$. On $\Pi$ one has an injection
$$
\Phi:\bigoplus^\ell \sO_{\Pi}(-1) \>>> \eta^*\sT,
$$
whose restriction to $\Pi_0$ is an isomorphism.
Writing $\eta_0:\Pi_0\to Y$ for the restriction of $\eta$ and
$$
\sO_{\Pi_0}(\mu)=\sO_{\Pi}(\mu)|_{\Pi_0},
$$
one considers the splitting 
$$
S^\nu\big(\bigoplus^\ell \sO_{\Pi_0}(-1)\big) \> S^\nu(\Phi) > \cong > \eta_0^*S^\nu(\sT)=\eta_0^* \sK_\nu\oplus
\eta_0^*\sT_\nu. 
$$
For $\ell_\nu=\rk(\sT_\nu)$ the projection to $\eta_0^*\sT_\nu$ defines a morphism
$$
\Psi: \Pi_0\>>> \Gr = {\rm Grass} (\ell_\nu , S^\nu (\C^\ell ))
$$ 
to the Grassmann variety $\Gr$ parameterizing $\ell_\nu$ dimensional quotient bundles. 
An ample invertible sheaf on $\Gr$ is given by the determinant of the 
universal quotient bundle, hence
$$
\sL_\nu=\eta_0^*\det(\sT_\nu)\otimes \sO_{\Pi_0}(\nu \cdot \ell_\nu)) 
$$
is the pullback of an ample sheaf on $\Gr$. Since  
$$
\sO_{\Pi_0}(\ell \cdot \nu \cdot \ell_\nu )=
\sO_{\Pi}(\nu\cdot \ell_\nu \cdot\Delta)|_{\Pi_0}=\sO_{\Pi_0},
$$
one finds $\sL_\nu^\ell=\eta_0^*\det(\sT_\nu)^\ell$. 
\begin{claim}\label{mu.4}
The image $\Psi(\eta_0^{-1}(y))$ is independent of $y$.
\end{claim}
\begin{proof}
Otherwise for some $\alpha>0$, divisible by $\ell$ there is a section of 
$\sL_\nu^{\alpha}$ which vanishes identically on $\Psi(\eta_0^{-1}(y))$, but not on $\Psi(\Pi_0)$. 
Then 
$$
\eta_0^*\det(\sT_\nu)^\alpha\otimes \sO_{\Pi_0}(-\eta_0^{-1}(y))
$$ 
has a non zero section. For some $\iota$ sufficiently large this section lies in
\begin{multline*}
H^0(\Pi,\eta^*\det(\sT_\nu)^\alpha \otimes \sO_{\Pi}(\iota\cdot\Delta-\eta_0^{-1}(y)))=\\
H^0(Y,\det(\sT_\nu)^\alpha\otimes S^{\iota\cdot\ell}({\sT^\vee}^{\oplus \ell})\otimes \sO_Y(-y)).
\end{multline*}
Since $\sT$ is poly-stable of degree zero and $\deg(\sT_\nu)=0$
there are no such sections.
\end{proof} 
To finish the proof of Claim \ref{mu.3} consider two points $y$ and $y'$ in $Y$.
By \ref{mu.4} 
$$
\Psi(\eta_0^{-1}(y))=\Psi(\eta_0^{-1}(y')),
$$
hence choosing suitable coordinates, the defining equations for $\pi^{-1}(y)\subset \BP^{\ell-1}$
and $\pi^{-1}(y')\subset \BP^{\ell-1}$ are the same.
\end{proof}
Let $\sB_\nu$ be the first sheaf in the Harder-Narasimhan filtration of 
$$
f_*\omega_Y(S)^\nu\otimes \sL^{-n\cdot\nu}
$$ 
and let 
$$
\mu(\sB_\nu)=\frac{\deg(\sB_\nu)}{\rk(\sB_\nu)} 
$$
denote the slope of $\sB_\nu$.
\begin{claim}\label{mu.5}
If $\mu(\sB_\nu)>0$, then (replacing again $Y$ by some covering, \'etale over $U$) 
for some $\mu \gg 1$ and for a general point $y\in Y$ there exists a section
of 
$$
f_*\omega_{X/Y}^{\mu\cdot\nu}\otimes \sL^{-n\cdot\mu\cdot\nu} \otimes \sO_Y(-\mu\cdot\nu\cdot y).
$$
\end{claim}
\begin{proof} 
Replacing $Y$ by a covering, we may assume that $d=\mu(\sB_\nu) > \nu$. 
Then the image of $\sB_\nu$ under the multiplication map
$$
S^\mu(\sT_\nu) \>>> \pi'_*\sO_{W}(\mu\cdot\nu).
$$
has slope larger than or equal to $\mu\cdot d$, hence the same holds true for $\sB_{\mu\cdot\nu}$. 
The Riemann-Roch Theorem for locally free sheaves on curves implies that
\begin{multline*}
\dim(H^0(Y, \sB_{\mu\cdot\nu}
\otimes \sO_Y(-\mu\cdot \nu\cdot y))) \geq \rk(\sB_{\mu\cdot\nu})
\cdot (\mu(\sB_{\mu\cdot\nu})-\mu\cdot\nu+1-g(Y))\\
\geq \rk(\sB_{\mu\cdot\nu})\cdot (\mu\cdot (d -\nu) + 1-g(Y)).
\end{multline*}
Obviously this is larger than one for $\mu \gg 1$.
\end{proof}
\begin{proof}[Proof of Proposition \ref{mu.1}]
By assumption $W\to Y$ is not birationally isotrivial, hence Claim \ref{mu.3} implies that
$\deg(\sT_\nu)>0$ for some $\nu>0$. Obviously this implies that $\deg(\sB_\nu)>0$ and
by Claim \ref{mu.5}
$$
f_*\omega_{X/Y}^{\mu\cdot\nu}\otimes \sL^{-n\cdot\mu\cdot\nu} \otimes \sO_Y(-\mu\cdot\nu\cdot y).
$$
has a non-trivial section. Since $\sL^2=\omega_S(Y)$ one obtains the condition ($*$).
\end{proof}
\begin{proof}[Proof of Proposition \ref{mu.2}]
As above we may assume that $\#S$ is even and write $F^{1,0}=\sL\otimes \sT$ with $\sT$ unitary,
and with $\det(\sT)=\sO_Y$. If $S\neq \emptyset$ \cite[Theorem 5.2]{VZ3} implies that 
$\varphi^*\sT$ is trivial for some covering $\varphi:\hat{Y}\to Y$. 
Hence we only have to consider the case where $Y=U$. In particular the genus of $Y$
is strictly larger than one.
 
The assumption $\ell >1$ implies that $\pi':W\to Y$ has one dimensional fibres. Let $\varsigma:\tilde W \to W$ 
be the normalization, and let $\tilde \pi : \tilde W \to Y$ the induced morphism. Furthermore we will write 
$\sO_{\tilde W}(\nu)=\varsigma^*\sO_{W}(\nu)$.

If $\deg(\sT_\nu)>0$ for some $\nu$, again $\deg(\sB_\nu)>0$ and Claim \ref{mu.5} implies the condition ($*$).
So by Claim \ref{mu.3} it remains to consider the case where all fibres of $\pi'$ are isomorphic,
hence the fibres of $\tilde \pi$ as well.\\

Assume first, that the genus of the fibres of $\tilde \pi$ is larger than or equal to $1$.
As well known (see \cite{Vie0}, for example), there exists an \'etale finite covering
$\hat{Y}\to Y$ such that $\tilde W \times_Y\hat{Y}$ is birational to $F\times \hat{Y}$. Since $g(\hat{Y}) > 1$, 
the surfaces $\tilde W \times_Y\hat{Y}$ and $F\times \hat{Y}$ are both minimal and of non negative Kodaira dimension,
hence they are isomorphic. 

Replacing $Y$ by $\hat{Y}$ we may assume that $\tilde W= F\times Y$. The image $F'$ of $F\times \{y\}$
in $\BP^{\ell-1}=\pi^{-1}(y)$ is independent of $y\in Y$, up to the action of $\BP{\rm Gl}(\ell,\C)$.
Since the automorphism group of $F'$ is finite, this implies that $\BP(\sT)=\BP^{\ell-1}\times Y$, and
$\sT$ is the direct sum of line bundles of degree zero.\\

It remains to consider the case that $\tilde \pi : \tilde W \to Y$ is a $\BP^1$ bundle, say
$\BP(\sE)$ for some locally free sheaf $\sE$. The invertible sheaf 
$\sO_{\tilde W}(\nu)$ has to be of the form $\sO_{\BP}(\sE)(r)\otimes {\pi'}^*\sN$ where
$\sO_{\BP}(1)$ is the tautological bundle and $\sN$ an invertible sheaf.
Replacing $Y$ by some \'etale covering, we may assume that $\sN$ is the $r$-th power of some invertible sheaf,
and changing $\sE$ we can as well assume that $\sN=\sO_Y$. Then we have inclusions
\begin{equation}\label{eqmu.1}
\sT_\nu \> \subset >> \pi'_*\sO_W(\nu) \> \subset >> \tilde \pi_* \sO_{\BP(\sE)}(r\cdot \nu)= S^{r\cdot \nu}(\sE)
\> \subset >> f_*\omega_Y(S)^\nu\otimes \sL^{-\nu}.
\end{equation}
The sheaf $\sT_\nu$ is a quotient of the poly-stable sheaf $S^\nu(\sT)$ of degree zero, hence poly-stable.

If the first step of the Harder-Narasimhan filtration $\sB_\nu$ of the right hand side
has a positive slope, Claim \ref{mu.5} implies the condition ($*$). 

Otherwise $\mu(\sB_\nu)=0$. If $\sE$ is not semi-stable, its Harder-Narasimhan filtration is of the form 
$$  
0\>>> \sM_1 \>>> \sE \>>> \sM_2 \>>> 0
$$
for line bundles $\sM_i$ with $\deg(\sM_1)> \deg(\sM_2)$, and $\deg(\sB_\nu)=0$ implies that
$0\geq \deg(\sM_1)$. Then the largest semi-stable
subsheaf of $S^{r\cdot \nu}(\sE)$ of degree zero is either zero or $\sM_1^{r\cdot \nu}$, a contradiction.

So $S^{r\cdot \nu}(\sE)$ is semi-stable of slope $r\cdot \nu \cdot \mu(\sE) \leq 0$. Since
$\sT_\nu$ is a subsheaf of $S^{r\cdot \nu}(\sE)$, this is only possible
if $\mu(\sE) = 0$.

Recall that the maximality of the Higgs field gives rise to a non trivial map
$$
S^2(\sT\otimes \sL) \> \subset >> S^2(f_*\omega_{X/Y}) \>>> \omega_Y(S).
$$
It factors through
$$
\sT_2\otimes \sL^2 \> \subset >> f_*\omega_{X/Y}^2 \>>> \omega_Y(S)=\sL^2
$$
and since $\sT_2$ is poly-stable of degree zero, the composite must split. So 
$\sT_2$ has a non trivial section. The inclusions in (\ref{eqmu.1}), for $\nu=2$, 
shows that $S^{2\cdot r}(\sE)$ has a section, which splits locally.
Hence $\sO_{\BP(\sE)}(2\cdot r)$ has a section whose zero divisor $D$ 
does not contain a fibre. If $D$ decomposes as
$$
D=\sum_{i=1}^s \alpha_i\cdot D_i,
$$
each $D_i$ is given by a section of $\sO_{\BP(\sE)}(\beta_i)\otimes \pi^*\sN_i$ for some invertible sheaf 
$\sN_i$ on $Y$, necessarily of non negative degree. Since
$$
\bigotimes_{i=1}^s \sN_i^{\alpha_i}=\sO_Y, 
$$
$\deg(\sN_i)=0$ and $D_{\rm red}$ is a section
of $\sO_{\BP}(\sE)(\beta)\otimes \pi^*\sN$, with $\deg(\sN)=0$. 

The canonical sheaf of $\BP(\sE)$ is given by 
$$
\omega_{\BP(\sE)/Y}= \tilde\pi^*\det(\sE)\otimes \sO_{\BP(\sE)}(-2)
$$ 
and one has an exact sequence
$$
0 \> >> \tilde{\pi}_* \omega_{\BP(\sE)/Y}(D_{\rm red}) \>>> \tilde\pi_*\omega_{D_{\rm red}/Y} \>>>
R^1\tilde\pi_*\omega_{\BP(\sE)/Y} \>>> 0.
$$
The right hand side is isomorphic to $\sO_Y$, and the left hand side
to
$$
\det(\sE)\otimes \sN\otimes \tilde\pi_* \sO_{\BP(\sE)}(\beta -2)=
\det(\sE)\otimes \sN\otimes S^{\beta-2}(\sE).
$$
So $\tilde\pi_*\omega_{D_{\rm red}/Y}$ is the extension of a semi-stable sheaf of degree zero
and of $\sO_Y$, hence again semi-stable of degree zero. This implies that $D_{\rm red}$ is \'etale over
$Y$.

Replacing $Y$ by an \'etale covering, we may assume that
$D$ is the sum of disjoint sections $D_1, \ldots , D_{2\cdot r}$, not necessarily distinct.

Then $\tilde\pi:\BP(\sE) \to Y$ has a section $s:Y\to \BP(\sE)$ with image $D_1$ and 
$$
s^* \sO_{\BP(\sE)}(-1)\cong  s^*\big(\omega_{\BP(\sE)/Y}(D_1)\otimes \tilde{\pi}^*(\sN_1^{-1}\otimes \det(\sE)^{-1})\big) \cong \sN_1^{-1}\otimes \det(\sE)^{-1}.
$$ 
This is only possible if the semi-stable sheaf $\sE$ is an extension 
$$
0\>>> \sN_1^{-1} \>>> \sE \>>> \sN_1 \otimes \det(\sE)\>>> 0 
$$ 
of invertible sheaves of degree zero. 
Then the graded sheaf for the Jordan-H\"older filtration of $S^r(\sE)$  
is the direct sum of invertible sheaves of degree zero. This property is inherited by any 
semi-stable subsheaf of degree zero, in particular using (\ref{eqmu.1}) for $\nu=1$, 
by $\sT$. Since $\sT$ is poly-stable, it must be a direct sum of invertible subsheaves of degree zero.
It remains to show, that over some \'etale covering of $Y$, the direct factors become isomorphic. 
\begin{claim}\label{mu.6}
If $\sT$ is the direct sum of invertible sheaves, then there exists an \'etale covering
$\varphi:\hat{Y}\to Y$ such that $\varphi^*\sT=\sO_{\hat{Y}}^{\oplus \ell}$.
\end{claim}
\noindent {\it Proof.} \,
By \cite[Lemma 3.2]{VZ2} one can assume that $\V$ and the decomposition $R^1f_*\C_V=\V\oplus \W$
are defined over some number field. Recall that $\L=\L^\vee$ and $\T=\T^\vee$, hence
one has isomorphisms of variations of Hodge structures 
$$
\END (\V)=(\L\otimes \T)^{\otimes 2}= \L^{\otimes 2}\otimes \T^{\otimes 2}=
(\END_0(\L)\oplus \C) \otimes \END(\T).
$$
As in the proof of \cite[Lemma 3.7]{VZ3}, the Lemma 3.5 in \cite{VZ3} implies that
the decomposition $\V=\L\otimes \T$ is defined over some Galois extension $K$ of $\Q$ with Galois group $G$.
Using \cite[Lemma 3.2]{VZ3} again, we may also assume that the decomposition
of $\T$ in rank one local sub systems is defined over the same field $K$.

Writing $\V_K$ and $\T_K$ for the corresponding $K$ local systems,
$\END(R^1f_*K_V)$ contains $\END(\V_{K})$, hence $\END(\T_{K})$.
Since $\T_K$ is the direct sum of rank $1$ local sub systems 
$\END(\T_{K})$ has the same property, as well as $\END(\T_{K})^\gamma$ for $\gamma\in G$. 

Consider the Weil restriction 
$$
\M_\Q=\sW_{K/\Q}(\END(\T_K))\subset \END(R^nf_*\Q_V).
$$
Since $\M_\Q\otimes K$ is the direct sum of rank one local systems, $\M_\Q$ is unitary.
Since it has a $\Z$-structure, \cite[Lemma 4.3]{VZ2} implies that it trivializes over a finite \'etale covering $\varphi: \hat{Y}\to Y$. Then $\varphi^*\END(\T)=\sO_{\hat{Y}}^{\oplus \ell^2}$, and replacing
$\hat{Y}$ by a degree two \'etale covering $\varphi^*\T=\sO_{\hat{Y}}^{\oplus \ell}$ .
\end{proof}
\section{The proof of Theorem \ref{in.7}}\label{pr}

In order to show that the condition ($*$) in Proposition \ref{mu.1} and the condition
$(**)$ both lead to contradictions we may replace $\hat{Y}$ by $Y$. So we will assume throughout this section that $f:X\to Y$ is a semistable family of manifolds and that at least one of the following
two assumptions holds true: 
\begin{enumerate}
\item[($*'$)]
For some $\nu$ sufficiently large and divisible by $2$
and for a point $y\in Y$ in general position there exists a section
of 
$$
f_*\omega_{X/Y}^\nu\otimes \omega_{Y}(S)^{-\frac{n\cdot\nu}{2}} \otimes \sO_{Y}(-\nu\cdot y),
$$
\item[($**'$)] $f:X\to Y$ is a family of curves and $F^{1,0}$ the direct sum
of at least two copies of an invertible sheaf $\sL$ with 
$\sL^2=\omega_{Y}(S)$. 
\end{enumerate}
For ($*'$) we will use methods from \cite{VZ1} which allow to
control the Kodaira-Spencer maps of the families. In particular we will have to
recall the main covering construction from \cite[Section 3]{VZ1}.
\begin{setup}\label{pr.1}
Let $\varphi:Y'\to Y$ be a finite covering and let $f':X'\to Y'$ the family obtained as pullback 
of $f:X\to Y$. Remark that the semi-stability of $f$ implies that
$X'$ is normal, with at most rational double points as singularities.
Consider a birational morphism $\delta:Z\to X'$, with $Z$ a manifold, and a finite 
Galois covering $W\to Z$ with Galois group $\Z/\nu$. So we have a diagram 
\begin{equation}\label{eqpr.1}
\begin{CD}
W \>\tau >> Z \> \delta >> X' \> \varphi' >> X\\
\V h VV \V g VV \V f' VV \V f VV\\
Y' \> = >> Y' \> = >> Y' \> \varphi >> Y.
\end{CD}
\end{equation}
We will write $\pi=\varphi'\circ\delta$. 
Let $\sM$ be an invertible sheaf on $Y'$, and let $\sigma$ be a section 
$$
\sigma \in H^0(Z,\delta^*\omega_{X'/Y'}^\nu\otimes g^*\sM^{-\nu}).
$$ 
We assume that:
\begin{enumerate}
\item[i.] $\tau:W\to Z$ is the finite covering obtained by taking the $\nu$-th root out of $\sigma$
(see \cite{EV}, for example).
\item[ii.] $g$ and $h$ are both smooth over $Y'\setminus T$
for a divisor $T$ on $Y'$ containing $\varphi^{-1}(S)$. 
Moreover $g$ is semi-stable and the local monodromy operators of $R^nh_*\C_{W\setminus h^{-1}(T)}$ 
in $t\in T$ are unipotent. 
\item[iii.] Let $\Delta'=g^*T$ and let $D$ be the zero divisor $\sigma$ on $Z$. Then
$\Delta'+D$ is a normal crossing divisor and $D_{\rm red}\to Y'$ is \'etale over $Y'\setminus T$.
\end{enumerate}
\end{setup}
$W$ might be singular, but the sheaf $\Omega_{W/Y'}^p(\log \tau^*\Delta')=\tau^*\Omega_{Z/Y'}^1(\log
\Delta')$ is locally free and compatible with desingularizations.
The Galois group $\Z/\nu$ acts on the direct image sheaves
$\tau_*\Omega_{W/Y'}^p(\log \tau^*\Delta')$. As in \cite{EV} or \cite[Section 3]{VZ1} one has the following description of the sheaf of eigenspaces.
\begin{lemma}\label{pr.2}  
Let $\Gamma'$ be the sum over all components of $D$, whose multiplicity
is not divisible by $\nu$. Then the sheaf
$$
\Omega^p_{Z/Y'}(\log (\Gamma'+\Delta'))\otimes \delta^*\omega_{X'/Y}^{-1}\otimes g^* \sM \otimes \sO_{Z}\big(\big[\frac{D}{\nu}\big]\big),
$$
is a direct factor of ${\tau}_*\Omega^p_{W/Y'}(\log {\tau}^*\Delta')$. Moreover the $\Z/\nu$ action on $W$
induces a $\Z/\nu$ action on 
of 
$$
\W=R^nh_*\C_{W\setminus \tau^{-1}\Delta'}
$$ 
and on its Higgs bundle. One has a decomposition of $\W$ in a direct sum of sub variations of Hodge structures, given by the eigenspaces for this action, and the Higgs bundle of one of them is of the form
$ G=\bigoplus_{q=0}^n G^{n-q,q}$ for
$$
G^{p,q}=R^qg_*\big(\Omega^{p}_{Z/Y'}(\log (\Gamma'+\Delta'))\otimes 
\delta^*\omega_{X'/Y}^{-1}\otimes \sO_{Z}\big(\big[\frac{D}{\nu}\big]\big)\big)\otimes \sM.
$$
The Higgs field $\theta_{p,q}:G^{p,q} \to G^{p-1,q+1}\otimes \Omega^1_{Y'}(\log T)$ is induced by the edge 
morphisms of the exact sequence
\begin{multline}\label{eqpr.2}
0\>>> 
\Omega^{p-1}_{Z/Y'}(\log (\Gamma'+\Delta'))\otimes {g}^* \Omega^1_{Y'}(\log T)\\
\>>> \Omega^{p}_{Z}(\log (\Gamma'+\Delta'))
\>>> \Omega^{p}_{Z/Y'}(\log (\Gamma'+\Delta')) \>>> 0,
\end{multline}
tensorized with $\delta^*\omega_{X'/Y'}^{-1}\otimes g^*\sM \otimes 
\sO_{Z}\big(\big[\frac{D}{\nu}\big]\big)$.
\end{lemma}
\begin{example}\label{pr.3}
Let us assume that the condition ($*'$) holds true, hence that there is a section $s$ of 
$$
\omega_{X/Y}^\nu\otimes f^*(\sL)^{-n\cdot \nu}\otimes \sO_Y(-\nu\cdot y)).
$$
Consider a desingularization $\hat{W}$ of the cyclic covering defined by $s$. Then $\hat{h}:\hat{W}\to Y$
will be smooth outside of a divisor $\hat{T}$, but not semi-stable. Choose $Y'$ to be a covering,
sufficiently ramified, such that the local monodromy operators of the pullback of $\hat{h}_*\C_{\hat{W}\setminus\hat{h}^*
\hat{T}}$ to $Y'$ are unipotent. Next choose $W'$ to be a $\Z/\nu$ equivariant desingularization of 
$\hat{W}\times_YY'$, and $Z$ to be a desingularization of the quotient.
Finally let $W$ be the normalization of $Z$ in the function field of $\hat{W}\times_YY'$.
So we constructed the diagram (\ref{eqpr.1}). 

For $\sM=\varphi^*\sL^{n}\otimes \sO_Y(y)$ consider the section $\varphi'^*(s)$
of 
$$
\varphi'^*\omega_{X/Y}^\nu\otimes f'^*\sM^{-\nu}=\omega_{X'/Y'}^\nu\otimes f'^*\sM^{-\nu}
$$
and the induced section $\sigma$ of $\delta^*\omega_{X'/Y'}^\nu\otimes g^*\sM^{-\nu}$.

The sum of the zero locus and the singular fibres will become a normal crossing divisor after
a further blowing up. Then one chooses $Y'$ larger, and one may assume that $Z\to Y'$ is semi-stable,
and that $Z$ and $D$ satisfy the assumption iii) in \ref{pr.1}. 
\end{example}
Remark that one has no controll on the critical locus $T$ of the family $Z\to Y'$.
Also, the sheaves $\Omega^p_{Z/Y'}(\log(\Gamma'+\Delta'))$ occurring in the description of
the Higgs bundles in \ref{pr.2} are not pullbacks of sheaves on $X$, and it is hard to describe
$\delta_*(\Omega^p_{Z/Y'}(\log(\Gamma'+\Delta')))$, even for $p=n$. However, if one forgets
the logarithmic poles along $\Gamma'$ one finds
$$
\delta_*(\Omega^n_{Z/Y'}(\log \Delta'))=\omega_{X'/Y'}={\varphi'}^*\omega_{X/Y}.
$$ 
As we will see in the proof of the next Proposition, this observation will allow
to construct certain logarithmic Higgs fields without poles in points of $T\setminus S$.

\begin{proposition}\label{pr.4}
Let $f:X\to Y$ be a semi-stable family of varieties over a curve $Y$, with a connected general fibre.
Let $\sL$ be an invertible sheaf with
$\sL^2=\omega_Y(S)$. Then for a general point $y\in Y$ and for all $\nu > 1$ 
$$
H^0(Y,f_*\omega_{X/Y}^\nu\otimes \sL^{-n\cdot\nu}\otimes \sO_Y(-\nu\cdot y))=0.
$$
\end{proposition}
\begin{proof}
If not choose $s$ to be the corresponding section of
$$
\omega_{X/Y}^\nu\otimes f^*(\sL^{-\nu}\otimes \sO_Y(-\nu\cdot y)),
$$
and perform the construction described in Example \ref{pr.3}. In particular for
$\sM=\varphi^*\sL^{n}\otimes \sO_Y(y)$ one obtains the diagram in \ref{pr.1}, satisfying the
properties i), ii), and iii), for the zero divisor $D$ of $\sigma$. So Lemma \ref{pr.2}
gives the description of the Higgs bundle of a particular sub variation of Hodge structures
of $R^nh_*\C_{W\setminus \tau^{-1}\Delta'}$. Using the notations introduced there, the sheaf
$$
G^{n,0}=g_*\big(\Omega^n_{Z/Y'}(\log (\Gamma'+\Delta'))\otimes 
\sO_{Z}\big(\big[\frac{D}{\nu}\big]\big)\otimes \delta^*\omega_{X'/Y'}^{-1}\big)\otimes \sM 
$$
contains
$$
\sH=g_*\big(\Omega^n_{Z/Y'}(\log \Delta')\otimes \otimes \delta^*\omega_{X'/Y'}^{-1}\big)\otimes \sM=
g_*(\omega_{Z/Y'} \otimes \delta^*\omega_{X'/Y'}^{-1})\otimes \sM.
$$
Since $f':X'\to Y'$ is a semi-stable family of $n$-folds over a curve, $X'$ has at most rational 
double points. Then 
$$
\delta_*(\omega_{Z/Y'} \otimes \delta^*\omega_{X'/Y'}^{-1})=\delta_*\omega_{Z/X'}= \sO_{X'},
$$ 
and
\begin{equation}\label{eqpr.4}
\sH=\sM=\varphi^*(\sL^n\otimes \sO_Y(y)).
\end{equation}
Let $(H=\bigoplus_{q=0}^nH^{n-q,q} ,\theta|_H)$ be the Higgs sub bundle of $(G=\bigoplus_{q=0}^nG^{n-q,q},\theta)$ generated by $\sH$. Then $H^{n,0}=\sH$ and 
$$
H^{n-q-1,q+1}={\rm Im}\big(\theta|_{H^{n-q,q}}:H^{n-q,q} \to G^{n-q+1,q+1}\otimes \omega_{Y'}(T)\big)\otimes \omega_{Y'}(T)^{-1}.
$$
In particular there is some $q_0\geq 0$ such that $H^{n-q,q}$ is an invertible sheaf for $q\leq q_0$ and
zero for $q>q_0$. 
\begin{claim}\label{au.6}
The image $\theta(H^{n-q,q})$ lies in $H^{n-q-1,q+1}\otimes \varphi^*\omega_Y(S)$, and
for $q\leq q_0$ 
$$
\deg(H^{n-q,q})=\deg(\varphi)\cdot \big((n-2\cdot q)\cdot\deg(\sL)+ 1\big).
$$
\end{claim}
\begin{proof}
Writing $\Delta=f^*(S)$ consider the tautological exact sequences
\begin{equation}\label{eqpr.3}
0\to
\Omega^{p-1}_{X/Y}(\log \Delta)\otimes {f}^* \Omega^1_{Y}(\log S)
\>>> \Omega^{p}_{X}(\log \Delta)
\>>> \Omega^{p}_{X/Y}(\log \Delta) \to 0,
\end{equation}
tensorized with $\omega_{X/Y}^{-1}=(\Omega^{n}_{X/Y}(\log \Delta))^{-1}$. Taking the edge morphisms
one obtains a Higgs bundle starting with the $(n,0)$ part $\sO_Y$. The sub 
Higgs bundle generated by $\sO_Y$ has $\omega_Y(S)^{-q}$ in degree $(n-q,q)$.
Tensorizing with $\sL^n\otimes \sO_Y( y)$ one obtains a Higgs bundle $H_0$ with
$$
H_0^{n-q,q}=\sL^n\otimes \sO_Y( y)\otimes \omega_Y(S)^{-q}=
\sL^{n-2\cdot q}\otimes \sO_Y( y),
$$
hence $\deg(H_0^{n-q,q})=(n-2\cdot q)\cdot\deg(\sL)+ 1$.

On the other hand, the pullback of the exact sequence (\ref{eqpr.3}) to $Z$ is a subsequence of
$$
0\to
\Omega^{p-1}_{Z/Y'}(\log \Delta')\otimes {g}^* \Omega^1_{Y'}(\log T)
\to \Omega^{p}_{Z}(\log \Delta')
\to \Omega^{p}_{Z/Y'}(\log \Delta') \to 0,
$$
hence of the sequence (\ref{eqpr.2}), as well. Then the Higgs field of $\varphi^*H_0$ with
$$
\theta': \varphi^* H_0^{n-q,q}\>\cong >> \varphi^* (H_0^{n-q-1,q+1}\otimes \omega_{Y}(S))
\> \subset >> \varphi^*H_0^{n-q-1,q+1}\otimes \omega_{Y'}(T)
$$
commutes with the edge morphism of the exact sequence (\ref{eqpr.2}), tensorized with
$$
\delta^*\varphi'^*\big(\omega_{X/Y}^{-1}\otimes f^* (\sL^n\otimes \sO_Y( y))\big)=
\delta^*\varphi'^*\Omega_{X/Y}^{-1}\otimes g^*\sM = \delta^*\Omega_{X'/Y'}^{-1}\otimes g^*\sM
$$
or with the larger invertible sheaf
$$
\delta^*\Omega_{X'/Y'}^{-1}\otimes g^*\sM\otimes \sO_Z\big(\big[\frac{D}{\nu}\big]\big).
$$
One obtains a morphism of Higgs bundles $\varphi^*H_0\to G$. By (\ref{eqpr.4}) and by the 
definition of $H$
$$
\varphi^*H_0^{n,0} = \varphi^*(\sL^n\otimes \sO_Y(y))= H^{n,0} \> \subset >> G^{n,0},
$$
so $H$ is the image of $\varphi^* H_0$ in $G$. Then $H^{n-q,q}\cong\varphi^*H_0^{n-q,q}$ for $q\leq q_0$
and $H^{n-q,q}=0$, otherwise.
\end{proof}
The Claim \ref{au.6} implies that the degree of $H$ is  
$$
\deg(\varphi)\cdot \big(q_0+1 + \big( (n+1)\cdot n - 2\cdot 
\sum_{q=0}^{q_0} q\big)\cdot\deg(\sL)\big),
$$
hence
\begin{equation}\label{eqpr.5}
\deg(H) \geq \deg(\varphi)\cdot\big(1 + \big( (n+1)\cdot n - (q_0+1)\cdot q_0 \big)\cdot\deg(\sL)\big)
> 0.
\end{equation}
By Simpson's  correspondence \cite{Sim} the Higgs bundle of a variation of Hodge structures
with unipotent local monodromy operators is polystable of degree zero. So
$H$ is a non-trivial sub Higgs bundle of a polystable Higgs bundle, contradicting
(\ref{eqpr.5}).  
\end{proof}
\begin{proof}[Proof of Theorem \ref{in.7}]
Assume that there exists a complex polarized sub variation of Hodge structures, 
with a maximal Griffiths-Yukawa coupling, satisfying the condition i) or ii) 
of Theorem \ref{in.7}. Lemma \ref{in.4} allows to choose such a $\V$ with a maximal Higgs field.
Writing $(F,\tau)$ for the Higgs bundle of $\V$, the assumption that $F^{n,0}$ is non isotrivial
implies by Proposition \ref{mu.1} that the condition ($*$) holds true over some \'etale covering 
of $Y$. Similarly, if $f:X\to Y$ is a family of curves, by \ref{mu.2} either ($*$) holds
true, or $F^{1,0}$ is a direct sum of several copies of a logarithmic theta characteristic.
In order to obtain a contradiction, we may replace $Y$ by this covering, hence assume that
either the condition ($*'$) or the condition ($**'$) holds for $f:X\to Y$. 

Proposition \ref{pr.4} tells us that there can not exist any family with a non-trivial
section of the sheaf
$$
f_*\omega_{X/Y}^\nu\otimes \sL^{-n\cdot\nu}\otimes \sO_Y(-\nu\cdot y)), 
$$
hence no family satisfying the condition ($*'$). 

It remains to exclude the case of a family of curves with ($**'$), hence of a family
whose variation of Hodge structures contains a sub variation with $(1,0)$ part
$$
F^{1,0}=\bigoplus^\epsilon \sL
$$
for some logarithmic theta characteristic $\sL$, and with $\epsilon \geq 2$.
Obviously each of the direct factors $\sL$ defines a sub Higgs bundle $\sL\oplus \sL^{-1}$
of degree zero, hence a rank two sub variation of Hodge structures $\V$ with a maximal Higgs field. 
As recalled in Theorem \ref{in.8}, by \cite[Theorem 2.12]{Moe} this forces $U$ to be a Teichm\"uller 
curve. For those \cite{McM1} (see also \cite[Lemma 3.1]{Moe}) excludes the existence of a second local sub system $\V'\neq \V$ in $R^1f_*\C_V$ with a maximal Higgs field, contradicting the assumption $\epsilon >1$.
\end{proof}

\end{document}